\documentclass[oneside,english]{amsart}
\usepackage[T1]{fontenc}
\usepackage[utf8]{luainputenc}
\usepackage{amsthm}
\usepackage{graphicx}

\makeatletter
\numberwithin{equation}{section}
\numberwithin{figure}{section}
\theoremstyle{plain}
\newtheorem{thm}{\protect\theoremname}
  \theoremstyle{plain}
  \newtheorem{cor}[thm]{\protect\corollaryname}
  \theoremstyle{remark}
  \newtheorem{rem}[thm]{\protect\remarkname}

\usepackage{amsfonts}\usepackage{amsthm}\usepackage{enumerate}\usepackage{mathrsfs}\usepackage{cancel}

\usepackage{MnSymbol}

\input xy\huge
\xyoption{all}

\DeclareFontFamily{OT1}{pzc}{}
\DeclareFontShape{OT1}{pzc}{m}{it}{<-> s * [1.10] pzcmi7t}{}
\DeclareMathAlphabet{\mathpzc}{OT1}{pzc}{m}{it}

\newcommand{\lc}{\mathcal{L}}

\usepackage{times}

\makeatother

\usepackage{babel}
  \providecommand{\corollaryname}{Corollary}
  \providecommand{\remarkname}{Remark}
\providecommand{\theoremname}{Theorem}

\begin{document}

\title{Taylor Expansion Proof of the Matrix Tree Theorem - Part I}

\author{Amitai Zernik}

\address{Einstein Institute of Mathematics,}

\address{The Hebrew University of Jerusalem }

\address{Jerusalem, 91904, Israel}

\address{email: amitai.zernik@gmail.com}
\begin{abstract}
The Matrix-Tree Theorem states that the number of spanning trees of
a graph is given by the absolute value of any cofactor of the Laplacian
matrix of the graph. We propose a very short proof of this result
which amounts to comparing Taylor expansions.
\end{abstract}
\maketitle

\section{\label{sec:Introduction}Introduction}

Call a matrix $L=(L_{ij})$ \emph{Laplace-like} if it is symmetric
and $\forall i\;\sum_{j}L_{ij}=0$. 

Let $M_{ij}(L)$ denote the submatrix of $L$ obtained by deleting
the $i$'th row and $j$'th column. The $i,j$ cofactor of $L$ is
$C_{ij}:=(-1)^{i+j}\det M_{ij}(L)$. It is not hard to see that if
$L$ is Laplace-like then all its cofactors are the same; we denote
their common value by \emph{$C(L)$.}
\begin{thm}
\label{thm:Kirchoff}If $L$ is an $n\times n$ Laplace-like matrix
then 
\begin{equation}
C(L)=(-1)^{n+1}\sum_{\tau}A_{\tau}(L)\label{eq:MTT}
\end{equation}
 where the sum is taken over all trees $\tau$ whose set of vertices
is $\{1,...,n\}$ and 
\begin{equation}
A_{\tau}(L)=\prod_{i<j,\{i,j\}\in\mbox{ Edges}(\tau)}L_{ij}\label{eq:Amplitude}
\end{equation}
.
\end{thm}
A special case of the above gives the well-known Matrix-Tree Theorem.
\begin{cor}
The number of spanning trees of an unoriented graph is equal to $|C(L)|$
where $L=A-D$ for $A$ the adjacency matrix of the graph and $D$
the diagonal matrix of valencies of the graph.
\end{cor}
The original formulation and proof of this beautiful theorem is due
to Kirchoff \cite{Kirchoff}. Since then, several different proofs
have been proposed. Our proof of Theorem \ref{thm:Kirchoff} is given
in Section \ref{sec:Proof-of-MTT} below, and is shorter than any
we've seen. Roughly speaking, it amounts to comparing the Taylor expansion
of both sides of eq (\ref{eq:MTT}). In the second part of this paper,
we shall use the same technique to prove a generalization of the Matrix
Tree Theorem for (not necessarily symmetric) matrices whose columns
sum to zero.

I would like to thank David Kazhdan, Nati Linial, Ori Parzanchevski,
Ron Peled, Ron Rosenthal and Ran Tessler for their suggestions and
comments. I am especially grateful to Nati Linial for suggesting reference
\textbf{\cite{combinatorial-ammtt}}.

\section{Proof of Matrix Tree Theorem\label{sec:Proof-of-MTT}}

We denote by $\lc$ the linear space of all Laplace-like matrices. 

Let $T(L)=\sum_{\tau}A_{\tau}(L)$. The tangent space to $\lc$ at
any point is spanned by $v_{ij}:=\frac{\partial}{\partial L_{ij}}+\frac{\partial}{\partial L_{ij}}-\frac{\partial}{\partial L_{ii}}-\frac{\partial}{\partial L_{jj}}$
for $i<j$. We will prove that $v_{ij}(C(L))=v_{ij}(T(L))$ for all
$i,j$ and $L$; the theorem then follows since $\lc$ is path connected
and $C(0)=T(0)=0$.

Let $L'=M_{jj}(L_{+})$, where $L_{+}$ is obtained from $L$ by adding
the $j$'th row to the $i$'the row and the $j$'th column to the
$i$'th column. $L'$ is an $(n-1)\times(n-1)$ Laplace-like matrix.
We will show that $v_{ij}(T(L))=T(L')$ and $v_{ij}(C(L))=-C(L')$.
The result then follows by induction on the size of the matrix $n$.

\begin{figure}
\includegraphics[scale=0.3]{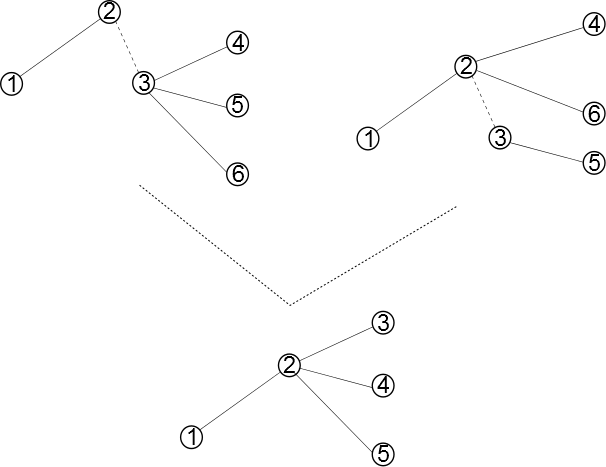}

\caption{\label{fig:contraction}Two distinct trees whose 2,3 contraction is
the same.}
\end{figure}

To see that $v_{ij}(T(L))=T(L')$, note that $v_{ij}(A_{\tau}(L))=\frac{\partial}{\partial L_{ij}}A_{\tau}(L)$
is nonzero only if $i,j$ are connected by an edge in $\tau$, in
which case we may contract the edge to produce a new spanning tree
$\tau'=\mbox{contract}_{ij}(\tau)$ with vertices $\{1,...,n-1\}$,
see Fig. \ref{fig:contraction}. Explicitly, we erase the vertex $j$
and the $i,j$ edge, and reconnect all the other neighbours of $j$
to the vertex $i$. We also relabel the vertices $j+1,...,n$ by $j,...,n-1$,
respectively. Fixing some $\tau'$, it is then easy to see that 
\[
\sum_{\{\tau|\mbox{contract}_{ij}(\tau)=\tau'\}}\frac{\partial}{\partial L_{ij}}A_{\tau}(L)=A_{\tau'}(L'),
\]
hence $v_{ij}(T(L))=T(L')$.

Next we show $-v_{ij}(C(L))=C(L')$. Indeed, 
\[
-v_{ij}(C(L))=-v_{ij}(C_{ii}(L))=\frac{\partial}{\partial L_{jj}}C_{ii}(L)=\frac{\partial}{\partial L_{jj}}C_{ii}(L_{+})=
\]

$=\frac{\partial}{\partial L_{jj}}\det(M_{ii}(L_{+}))=\det(M_{j-1,j-1}(M_{ii}(L_{+})))=\det M_{ii}(M_{jj}(L_{+}))=$

\[
=\det M_{ii}(L')=C(L')
\]
which completes the proof.\qed
\begin{rem}
Note even if $L$ is not symmetric (but its rows and columns sum to
zero) it is still true that $-v_{ij}(C(L))=C(L')$. This has led us
to suspect that a non-symmetric version of the theorem should exist.
In the second part of this paper we discuss a more general, non-symmetric
formulation of the Matrix Tree Theorem that we managed to prove using
a very similar technique.
\end{rem}
\bibliographystyle{plain}
\bibliography{matrix_tree_bibilography}

\end{document}